\documentclass[11pt]{article}
\usepackage{graphicx}
\usepackage{amssymb}

\usepackage[
            pdftex,
           colorlinks=true,
            urlcolor=blue,       
            filecolor=blue,     
            linkcolor=blue,       
            citecolor=blue,         
            pdftitle= {Title},
            pdfauthor={Author},
            pdfsubject={Subject},
            pdfkeywords={Key}
            pagebackref,
            bookmarksopen=true]
            {hyperref}
\usepackage{hyperref}

\usepackage[T2A,T1]{fontenc}
\newcommand{\Sh}{\mbox{\usefont{T2A}{\rmdefault}{m}{n}\CYRSH}}

\newcommand{\sh}{\textrm{sh}}
\newcommand{\Alb}{\textrm{Alb}}
\newcommand{\Aut}{\textrm{Aut}}
\date{}
\title{Igor Rostislavovich Shafarevich: in Memoriam}
\author{Igor V. Dolgachev}
\begin{document}
\maketitle
The prominent Russian mathematician Igor Rostislavovich Shafarevich passed away on February 19, 2017. 
He has made an outstanding contribution to number theory, algebra, and algebraic geometry. 
The influence of his work on the development of these fields in the second half of the 20th century 
is hard to overestimate. Besides the fundamental results authored by him and his 
collaborators, he single-handedly created a school of Russian algebraic geometers and number theorists, 
many of his numerous students consider their time spent under his guidance as the happiest time in their 
life as mathematicians. Shafarevich was awarded the Lenin prize in 1959 for his work on 
the inverse Galois problem, he was elected to the Russian Academy of Sciences in 1958 as a 
correspondent member and as a full member in 1991. Shafarevich was also a foreign member of 
the Italian Academy dei Lincei, the  German Academy Leopoldina, the National Science Academy of 
the USA (from which he resigned in 2003 as a protest against the Iraq War), 
a member of the London Royal Society, and received a honorary doctorate from the University of
 Paris. Shafarevich was an invited speaker at the International Congresses of 
 Mathematicians in Stockholm (1962) and Nice (1970). His name is associated with such 
 fundamental concepts and results in mathematics as the Shafarevich-Tate group, 
 Ogg-Shafarevich theory,  Shafarevich map, Golod-Shafarevich Theorem, 
 Golod-Shafarevich groups and algebras, Deuring-Shafarevich formula, and several Shafarevich Conjectures. 
 His influential textbooks in algebraic geometry and number theory (jointly with Zinovy Borevich) 
 have been translated into English and served as an introduction to these subjects for 
 several generations of mathematicians. His book, "Basic Notions in Algebra" \cite{ShafarevichBasic}, provides a bird's-eye view of algebra, 
 revealing its vast connections with many other fields of mathematics, and has 
 become a favorite book in the subject for many mathematicians. 
 I quote from the preface to a collection of papers ``Arithmetic and Geometry''
  published in two volumes by Birkh\"auser in 1983 and edited by M. Artin 
  and J. Tate \cite{Arithmetic}: `Igor Rostislavovich Shafarevich has made 
  outstanding contributions in number theory, algebra, and algebraic geometry. 
  The flourishing of these fields in Moscow since World War II owes much to his influence. 
  We hope these papers, collected for his sixtieth birthday, will indicate to him the great 
  respect and admiration which mathematicians throughout the world have for him.'' 

In the preface to \cite{ShafCW1}, Shafarevich writes, ``At the end of the 
sixties the perception of life began to change. The passiveness of thinking and muteness 
became felt as irresponsibility. This new feeling seemed to turn me onto another road. 
Otherwise, I would stay till the end of my life in my profession as a mathematician, 
and my interest in history would remain as a hobby. Instead of this, I had acquired 
the second working profession to which I devoted with more and more strength." 
The subsequent non-mathematical activity that led to his numerous publications   
on social issues had at the same time tarnished and magnified his reputation among 
different layers of society in Russia and the West. 

\subsection*{Biography} Igor Rostislavovich Shafarevich was born in 1923 in the 
Ukrainian town Zhitomir. The name of the town is explained by the old Russian word "zhito", which means ``rye''. The same town was the birth-town for many famous Russians, for example, the pianist Svyatoslav Richter who remained  a life-long friend of Shafarevich. 

Shafarevich's father,  Rostislav Stepanovich graduated from the mathematical 
department of Moscow State University (MGU) and, after moving to Moscow, lectured in 
theoretical mechanics at one of the Institutes of Higher Learning. 
His mother Julia Yakovlevna was a philologist and a gifted pianist. 
Apparently, she shared with his son her lifelong passion for classical music and Russian literature. 
Igor's first serious interest as a child was in history, to which he was devoted 
till the end of his life. His other love was mathematics. Still at school, he took exams 
in mathematics at MGU from which he had graduated in 1940 at the age 17. 
Although he did not have a formal thesis adviser, his advisor for the master thesis was 
Boris Nikolaevich Delone. Other mathematicians whom he acknowledged as his mentors were 
Israel Moiseevich Gelfand and Alexander Gennadievich Kurosh. He had finished graduate school at MGU with 
a Ph.D. dissertation `On normiering of topological fields' in 1943 at age 20. 
During World War II, along with some of the university's faculty, he was evacuated to Ashkhabad and later to Kazan. 
After returning to Moscow, he defended his second thesis (a Russian version of German Habilitation) 
in 1946. In his thesis, he described all p-extensions of the field of p-adic numbers 
and non-ramified extensions of the fields of algebraic numbers. His doctoral committee 
included such prominent Russian mathematicians as Dmitry Konstantinovich Faddeev, 
Anatoly Ivanovich Maltsev, and Nikolai Grigorievich Chebotarev. 
After the defense of his thesis and until his death, he was a member of the Steklov 
Institute of Mathematics. Also, since 1944, he was teaching at MGU, where in 
the sixties he founded his famous seminar in Algebraic Geometry. In 1975, he was dismissed from the university due to his dissident activities. 
His seminar had been moved to the Steklov Institute, where it still meets on Tuesdays. 
For many years, Shafarevich was directing the Algebra section of the Institute and 
was credited to the worldwide renowned center of mathematical activity in algebra, algebraic geometry and number theory. Although he was sometimes addressed by his students as a ``boss'', there was never anything bossy in his relationship with his students, colleagues and ordinary Russian people who later were coming  to him for an advice on social issues. He always  respected   his numerous students and colleagues, treated  them as equal, and was  ready to help them in their mathematical careers and difficult periods of their life. Some of them were his true friends with whom he shared his passion for mountains hikes  and who helped him in his dissident activity. 

Shafarevich's scientific honesty is clearly revealed in his mathematical writings. 
His attribution of known results and historical references should serve as 
instructive examples for mathematicians of later generations. On several occasions, 
he stood up to express critical opposition to the weak or erroneous theses in the mathematics department at MGU 
(including the Habilitation thesis of his former student A. Zhizhenko, now a full member of the Russian Academy of Science, who became a Soviet bureaucrat).

\subsection*{Students} Since late forties, Shafarevich began advising Ph.D. dissertations. 
If he were not dismissed from the University, his list of students would be much larger. 
The following is, hopefully, a complete list of his Ph.D. students. 
Together with the descendants, the list contains more than 300 names.

\begin{table}
\centering
\caption{Ph.D. Students}
\begin{tabular}{llll}\\
Abrashkin& Victor& MGU&1976\\
Arakelov& Suren J.&MGU& 1974\\
Averbuch& Boris G.&MGU& 1964\\
Belyi& Gennady V.& MIAN&1979\\
Berman& Samuil D. &MGU&1952	\\
Demyanov& V. V. &MGU&1952 \\
Demushkin& Sergei, &MGU&1959\\
Dolgachev& Igor V. &MGU&1970\\
Drozd& Yurii A &MGU&1970\\ 
Gizatullin& Marat H. &MGU& 1970\\
Golod& Evgeny S. &MGU&1960\\
Koch& Helmut &MGU& 1964\\
Kolyvagin& Victor A. &MGU&1981\\
Kostrikin& Alexsei I. &MIAN&1960\\ 	
Kulikov& Valentine S.&MGU& 1975\\	
Kulikov& Viktor S. &MGU&1977\\
Nikulin& Vyacheslav V. &MGU&1977\\
Lapin&	Andrei I. &MGU&1952\\
Manin& Yuri I. &MGU&1961\\
Markshaitis& Gamlet N.&MGU&?\\ 
Medvedev& P. A.&MGU&?\\
Milner& A.A.&MGU&?\\
Neumann& Olaf&MGU&1966\\
Pavlov& ? &MGU&?\\
Parshin& Alexei N. &MGU&1967\\
Rudakov& Alexei N.&MGU&?\\
Shabat& George B.&MGU& 1976\\
Todorov& Andrey N.&MGU&1976\\
Tyurina& Galina N.&MGU& 1963 \\
Tyurin& Andrei N. &MGU&1965\\
Vvedenskii& Oleg N.&MGU&1963\\ 
Zhizhchenko& Alexei B.&MGU&1958
\end{tabular}
\end{table}

\subsection*{Scientific work: Number Fields} In his Habilitation dissertation 
Shafarevich studied non-abelian $p$-extensions of local and global fields. 
For example, he proved that given a finite degree $n$ extension of the field $\mathbb{Q}_p$ of 
rational $p$-adic numbers that does not contain $p$-roots of unity, the Galois group of 
its finite $p$-extension is a quotient of a free group with $n+1$ generators \cite{Shafarevich1}. 
For this work, Shafarevich was awarded the prize of the Moscow Mathematical Society. 
In his next work, he made a major contribution to number theory by giving an explicit 
formula for the local symbol $(\frac{\alpha,\beta}{p})$ \cite{ShafarevichReciprocity}. 
The formula is reminiscent of a familiar formula for the residue of a differential on a Riemann surface. 
The theory developed in his dissertation gave a new approach to the global and local 
class theory (see \cite{Koch}). His next work was even more impressive. 
In paper \cite{ShafInverse} of 1954, Shafarevich solves the inverse Galois problem 
for solvable groups in the case of fields of algebraic numbers. 
A gap in the proof of this fundamental result, pointed out much later by H. Koch and A. Schmidt 
was fixed by Shafarevich in 1980 in one of the footnotes to his 
Collected Works  \cite{ShafarevichCP}, p. 752.  
The proof was based on his earlier paper on the construction of $p$-extensions of algebraic number 
fields and uses new pioneering methods of homological algebra developed around this 
time by D. K. Faddeev. A complete proof using new tools can be found in the book \cite{Neukirch}.

The next problem addressed by Shafarevich was the problem of embedding of local and global fields $k$. 
Given a Galois extension $L/k$ with Galois group $G$ and its Galois subextension $K/k$,
 the subgroup of $G$ fixing elements from $K$ is a normal subgroup of $G$, with quotient group $G'$ 
 isomorphic to the Galois group of $K/k$. Given a surjective homomorphism of groups $G\to G'$, 
 the embedding problem asks whether there exists an embedding of a Galois extension $K/k$ 
 with Galois group $G'$ into a Galois extension $L/k$ with Galois group $G$ that realizes 
 the surjective homomorphism as the quotient map. In the case when $G$ is abelian and 
 the surjection $G\to G'$ with kernel $H$ makes $G$ a semi-direct problem 
 $H\rtimes G'$, the problem was solved in 1929 by A. Scholz. Shafarevich generalizes 
 this result to the case where $H$ is a nilpotent group of a certain class. 
 He returns to the embedding problem later in a joint work with his former 
 student Sergei Demushkin, first considering the case of local fields \cite{ShafDem59} and, 
 later, the case of global fields \cite{ShafDem62}. 

In 1963, Shafarevich published an important paper 
in Publicationes Mathematique IHES \cite{ShafIHES} 
(a rather rare event after World War II when  a Soviet mathematician  publishes in a Western journal) 
on the problem of $p$-extensions of algebraic number fields by considering 
finite extensions of these fields with a fixed set $S$ of ramified divisors. 
In the case when $S$ is the empty set, Shafaevich shows that the minimal number 
$d$ of generators of the Galois group of the extension and the number $r$ of minimal 
relations between generators satisfies inequality $r\le d+\rho$, where $\rho$ is the 
number of generators of the group of units of the field. 

At his talk at the ICM in Stockholm, he remarks that if one proves that $r(G)-d(G)\to \infty$, where the limit is taken over the set of all $p$-groups, then the class field tower problem on the existence of infinite unramified extensions of an algebraic number field has a negative solution. In a joint work with his former student  Evgeny  Golod \cite{Golod1} he proves that the limit is in fact goes to infinity solving in this way a classical fundamental problem  in number theory of more than 40 years old.  

\subsection*{Scientific work: Elliptic curves}  The transition of Shafarevich's interests 
from number theory to algebraic geometry was rather smooth and was based on his, now-famous 
work on elliptic curves. Already in 1956, in his talk at the Third Congress of 
Soviet Mathematicians, he pointed out the analogy between the problem of embedding of 
algebraic number fields and the problem of classification of elliptic curves over such fields. 
Both problems employ the local-to-global approach: find a solution for all completions of the field and determine whether it yields a 
solution over a global field.  In the case of elliptic curves, this leads 
to a question of whether the set of elliptic curves with a fixed absolute invariant 
isomorphic to a fixed curve over all completions of the field is finite. 
In a short announcement note \cite{ShafElliptic1} published in Doklady AN SSSR, 
he shows that the set of elliptic curves isomorphic to a fixed curve over some extension 
of the ground field forms a group that admits a cohomological interpretation as 
the first Galois cohomology group with coefficients in the group of points of the Jacobian curve. 
The fact that such a set forms an abelian group was not new; in the case when the ground field is 
the field of real numbers, it was discovered by Francois Ch\^{a}telet in 1947, 
whose construction uses the same cocycles. In 1955, A. Weil extended this result to the case of abelian varieties of arbitrary dimension, 
although he did not give a cohomological interpretation of the group. 
The paper of S. Lang and J. Tate of 1958 gives a foundation of the theory of principal 
homogeneous spaces over an abelian variety based on 
its cohomological interpretation (without reference to Shafarevich's paper). 
They call the group the Ch\^{a}telet group, later known under the name the Weil-Ch\^{a}telet group. 
In the same paper, Shafarevich proves that the subgroup of the Weil-Ch\^{a}telet 
group of elements that admit a point of degree $n$ over the ground field (hence. birationally isomorphic 
to an elliptic curve of degree $n+1$ in $\mathbb{P}^n$) and isomorphic to its Jacobian curve over 
all completions of the field is a finite group. In the subsequent paper \cite{ShafExponents} in Doklady, 
Shafarevich proved the existence of elliptic curves of arbitrary degree $n$ not isomorphic 
to any curve of smaller degree, giving a solution to an old problem in the theory of 
diophantine equations.

In 1967, during his stay in Paris, Shafarevich cooperated with John Tate 
to construct examples of elliptic curves over the functional field $k(t)$ with 
finite field $k$ whose Mordell-Weil group of rational points has arbitrarily 
large rank \cite{ShafTate}. The analogous statement, where the field $k(t)$ is 
replaced with the field $\mathbb{Q}$ of rational numbers, is still a conjecture in number theory.

In 1961, Shafarevich published a paper devoted to a systematic study of 
the Weil-Ch\^{a}telet group $H^1(K,A)$ of an abelian variety $A$ over a field $K$ of algebraic 
functions in one variable over an algebraically closed field $k$. 
Thus, he divides this study into three parts by determining the structure of three groups: 
the local group $H^1(K_\mathfrak{p}, A_{\mathfrak{p}})$, where the field $K$ is 
replaced by its completion, the kernel and the cokernel of the restriction homomorphisms 
to the product of these groups with respect to the set of all completions of $K$. 
The kernel group (where $K$ is replaced by an arbitrary field) was later named  
the Tate-Shafarevich group (the order is taken according to the Cyrillic alphabet). 
Shafarevich's contribution to the theory of elliptic curves 
was specially honored by a common acceptance of using the Cyrillic letter for its notation $\Sh(A)$. 
A similar theory, and about the same time, was independently developed by 
Andrew Ogg in Berkeley.  Later on, the theory was given a more modern approach by 
Grothendieck who gave a cohomological interpretation of $\Sh(A)$ as the 
first \`etale cohomology group with coefficients in a sheaf over a curve $C$ with the 
field $k(C)$ of rational functions isomorphic to $K$, which is represented by the N\'eron model of $A$. 
The main result of Ogg and Shafarevich on the structure of $\Sh(A)$ is now known as the 
Grothendieck-Ogg-Shafarevich formula. In their work, Ogg and Shafarevich restricted 
themselves only to the part of the Weil-Ch\^{a}telet group prime to the characteristic of $k$. 
The subsequent work of several people, including Oleg Vvedensky, a former student of Shafarevich, 
finished the work by settling the $p$-part \cite{Berta}. 

It is remarkable that the last published work of Shafarevich, when he was 
90 years old, was in number theory. In \cite{ShafDiscriminant}, he gives a 
new proof using the theory of modular forms of Stark's theorem that there are 
only nine imaginary quadratic fields with class number one.

\subsection*{Scientific Work: Algebraic Geometry} Shafarevich was always interested in algebraic geometry. 
For example, in 1950, he authored an article on algebraic geometry in the Russian Encyclopedia. 
In his paper \cite{ShafExponents}, he referred to a paper by F. Enriques of 
1899, which contains some geometric analogs of some of his results. 
 It should be noted that algebraic geometry and the theory of algebraic functions in one variable 
were always outside the interests of Russian schools in mathematics. 
 The only textbook on this topic was Chebotarev's book \cite{Chebotarev}, published in 1948, which 
 gives an exposition of the algebraic theory of algebraic curves. In 1961-1963, 
 Shafarevich and a group of his students ran a seminar on algebraic surfaces 
 whose goal was to revive some of the classical works of Italian algebraic 
 geometers from a modern point of view. The new techniques based 
 on topological methods and the use of the new theory of cohomology 
 of algebraic coherent sheaves developed earlier by Jean-Pierre Serre were 
 common tools in their work. The same activity was also undertaken at about the same 
 time by Oscar Zariski and David Mumford at Harvard and Kunihiko Kodaira at Princeton. 
 A book 'Algebraic surfaces' had appeared in Russian in 1965 and had been translated into 
 English in the same year. For many years, this book has been the primary source 
 for learning the classification of algebraic surfaces from a modern point of view. 
 Shafarevich himself contributed two chapters to the book. 
 In one of them, he translated his previous work on principal homogeneous spaces of 
 elliptic curves into geometric language, in particular, reconstructing 
 Enriques' work on elliptic surfaces. In another chapter, he gave a modern proof 
 of Enriques's criterion of ruledness of algebraic surfaces. As his students acknowledge, 
 his influence on the book as a whole was much greater than just contributing two chapters. 
 A very appropriate epigraph chosen for the book reflects very well 
 Shafarevich's admiration of classical works ``Aischylos said that his tragedies were 
 leftovers from great feasts of Homer.''

In I971, Shafarevich turned his attention to the study of complex K3 surfaces, 
which represent the most interesting two-dimensional analogs of elliptic curves. 
Their occurrence in many areas of mathematics and even mathematical physics is really remarkable. 
K3 surfaces share one common property with elliptic curves: the existence of a unique, 
up to proportionality, holomorphic differential form of highest degree. However, they 
differ from elliptic curves by the property that they are simply connected. 
It is a simple fact that the complex structure of an elliptic curve is determined by its 
periods, i.e., the values of integrals of its holomorphic form on a basis of 1-homology of the curve. 
Considered as a vector 
$(\int_{\gamma_1}\omega,\int_{\gamma_2}\omega)$ modulo proportionality and modulo of the group 
$\textrm{SL}_2(\mathbb{Z})$ acting via basis changes, it is a point 
in $\mathbb{C}$ that determines the curve up to isomorphism. 
The proof of this fact follows easily from representing an elliptic curve as 
the quotient of $\mathbb{C}$ by the lattice spanned by the periods. 
The absence of this representation for K3 surfaces made  Andr\'e Weil's guess that the 
periods of K3 surfaces should also determine their holomorphic structure 
seemed to be too daring to attempt to prove. Weil himself recognized this 
by K3 surfaces:

 ``il s'agit des vari\'etes k\"ahleriennes dites K3, ainsi nomin\'ees en l'honneur de 
 Kummer, K\"ahler, Kodaira et de la belle montagne K2 au Cachemire.'
 
  Nevertheless, the joint work of Shafarevich and Ilya Iosephovich Pyatetsky-Shapiro 
  has done exactly this. They proved that a projective complex algebraic K3 surface 
  is uniquely determined by its vector of periods modulo proportionality and changes of 
  a basis in the subgroup of the 2-homology group orthogonal to the class of its hyperplane section. 
  This result became known as the \emph{Global Torelli Theorem} for algebraic K3 surfaces 
  named after an Italian algebraic geometer, Ruggiero Torelli, who proved a similar result 
  for algebraic curves \cite{ShafPyatetsky1}. A corollary of this theorem allowed 
  them to reduce the study of the automorphism group of a K3 surface to the study of 
  some arithmetical property of an integral quadratic intersection
   form of algebraic cycles on the surface. This became an essential tool 
   in subsequent and continuing extensive study of automorphism groups of K3 surfaces. 
  
  The absence of topological and analytical methods for studying K3 surfaces 
  defined over fields of positive characteristic seemed to be an insurmountable obstacle 
  for extending the study of K3 surfaces in this case. 
  A paper by Michael Artin \cite{Artin} (which Shafarevich 
  acknowledged to me to be one of the most beautiful papers he had read in his life) was 
  a breakthrough in this direction. In this paper, Artin introduced the periods of supersingular K3 
  surfaces, the surfaces that are distinguished by the property that they have the maximum possible 
  number of linearly independent algebraic cycles. 
  In a long series of influential papers with his former student Alexei Rudakov, 
  Shafarevich undertook a comprehensive study of K3 surfaces over fields of positive characteristic. 
 For example, they prove the unirationality of supersingular K3 surfaces over a field of 
 characteristic two, prove non-degeneracy of supersingular K3 surfaces, 
 the absence of non-trivial regular vector fields on K3 surfaces, and lay the foundations for the 
 theory of inseparable morphisms of algebraic varieties. 
Using the non-degeneracy results of Shafarevich and Rudakov, 
Arthur Ogus was able to prove a Global Torelli Theorem for supersingular K3 surfaces over 
fields of odd characteristic. 

The Global Torelli Theorem for K3 surfaces, together with the surjectivity of the period map
 for complex algebraic
K3 surfaces, as proved by his former student Andrei Todorov, allow one to construct a coarse 
moduli space for algebraic K3 surfaces as an arithmetic quotient of a Hermitian symmetric 
domain of orthogonal type. Apparently, Shafarevich was interested in the theory of 
arithmetic groups and automorphic functions for a long time. In 1954, 
he wrote a preface and edited the Russian translation of Siegel's book \cite{Siegel}. 
In his paper with Pyatetski-Shapiro \cite{ShafPyatetsky2}, he studies a pro-algebraic variety 
with the field of rational functions equal to the limit of the fields of automorphic functions 
of subgroups of finite index of a discrete arithmetic group of automorphisms of a bounded symmetric domain. 
The second volume of his `Basic Algebraic Geometry' ends with a discussion of a problem of 
uniformization of algebraic varieties and makes his famous Shafarevich Conjecture 
that suggests that the universal cover of a complex projective variety $X$ must be holomorphically convex.
 In other words, Shafarevich conjectured that the universal cover admits a proper map to a Stein manifold 
 with connected fibers. In another reformulation, due to 
 Janos Koll\'ar, there must be a proper map $\sh_X:X\to \Sh(X)$ onto a normal variety $\Sh(X)$ 
 with connected fibers that contracts all closed subvarieties $Y$ of $X$ such that 
 the natural homomorphism of the fundamental group $\pi_1(Y')$ of a a resolution of singularities of $Y$ to the fundamental group $\pi_1(X)$ has finite image.
 Koll\'ar named a map with this property the Shafarevich map. 
 Koll\'ar's monograph \cite{Kollar} contains an extensive study of 
 the Shafarevich Conjecture and culminates with a proof of the existence of a birational map $\sh_X'$ 
 with similar properties. The Shafarevich conjecture is closely related to the 
 group-theoretical properties of the fundamental group $\pi_1(X)$, for example, 
 the existence of its faithful representation in a simple compact Lie group with dense image.  
 
 The Shafarevich map $\sh_X$ should be considered as a non-abelian generalization of the 
 Albanese map $a_X:X\to \Alb(X)$ that has the same property with respect to abelian unramified 
 covers of $X$. In his popular article in Mathematical Intelligencer in 
 2009 \cite{ShafNonabelian}, Shafarevich proposed that the deepest challenges 
 of modern mathematics can be summed up as a ``non-abelianization of mathematics. 
 He acknowledged that the ``non-abelian mathematics of the 
future'' philosophy also inspired him when he started his work in mathematics.

The combined interest of Shafarevich in number theory and algebraic geometry 
is explained by many close analogies between the two theories that go back to Leopold 
Kronecker and David  Hilbert. Shafarevich's talk at the International Congress of Mathematicians 
in Stockholm in 1962 is entirely devoted to the connections between the two fields. In particular, he stated two very influential conjectures in his talk. The analog of the Hermite conjecture about the finiteness of the number of finite extensions of an algebraic number field with the fixed discriminant becomes his  conjecture about the finiteness of the set of algebraic curves of fixed genus $g > 0$ over a number field $k$ with fixed discriminant and an analog of Minkowski's theorem that there are no unramified extensions of $\mathbb{Q}$  that now states that there are no smooth families of curves of positive genus over $\textrm{Spec}(\mathbb{Z})$. The attempts to prove the firt conjecture played a crucial role in Falting's proof of the Mordell Conjecture.

The beginning of the sixties was a time when many algebraic geometers of the 
present and earlier generations had to reeducate themselves in learning the 
new language of algebraic geometry was developed by the fundamental work of 
Alexander Grothendieck. Bombay Lectures of Shafarevich on minimal models of 
two-dimensional schemes over a discrete valuation ring \cite{ShafBombay}, 
together with Mumford's Lectures on curves on algebraic surfaces \cite{Mumford} 
were instrumental tools for accomplishing this goal.

In \cite{ShafKazan}, Shafarevich stated a conjecture: the set of 
Picard lattices of K3 surfaces defined over a field of algebraic numbers of 
degree $n$ over $\mathbb{Q}$ is a finite set. He proved this conjecture 
for K3 surfaces with maximal Picard number equal to 20. 
He also proved a geometric analog of this conjecture for one-dimensional 
families of Kummer surfaces. In a paper \cite{Shaf1996} published in the same year, 
he studies the Shimura variety of abelian surfaces with quaternionic multiplication 
(fake elliptic curves) and proves that the number of isomorphism classes of non-constant 
fake elliptic curves defines over an extension $K/\mathbb{C}(t)$ of degree $\le n$ is finite.

\subsection*{Scientific Work: Algebra} The work of Shafarevich in number theory led him to 
some fundamental problems in group theory. For instance, the solution to the problem of the 
existence of an infinite tower of class field extensions led him and Evgeny Solomonovich Golod to 
proving that  $r > (\frac{d-1}{2})^2$, where $r$ is the smallest number of 
generators of a $p$-group $G$ and $d$ is the smallest number of its generators. 
It is known that the numbers $r$ and $t = r-d$ 
can be interpreted in terms of the group cohomology as 
$r= \dim H^1(G,\mathbb{Z}/p\mathbb{Z})$ and $t = \dim H^2(G,\mathbb{Z}/p\mathbb{Z})$. 
Thus, the Golod-Shafarevich inequality becomes an equality on the Betti numbers $b_i$ 
of the graded algebra of cohomology $H^*(G,\mathbb{Z}/p\mathbb{Z})$. 
The main implication of the Golod-Shafarevich inequality (later improved by E. Vinberg and P. Roquette 
to the form $r\le d^2/4$) is that the small number of relations compared to the 
number of generators implies that the group must be infinite. In this way, an analogous statement  
toin different categories can be proved by similar methods and is referred to as the 
Golod-Shafarevich Theorem. This also led to the definition of the \emph{Golod-Shafarevich group}
 as a $p$-group with certain properties of its presentation, which implies that the group is infinite. 
 There has been an extensive study of Golod-Shafarevich groups and their analogues in other categories. Also there are new applications of the Golod-Shafarevich theory. For example,  Alexander Lubotsky proved that the fundamental group of a hyperbolic 3-manifold of finite volume contains a Golod-Shafarevich subgroup of finite index.

In 1964-66, Shafarevich ran a seminar at the Steklov Institute on 
Cartan's classification of simple transitive transformation Lie pseudogroups. 
A result of this seminar is a joint paper of Shafarevich and his former student, Alexei Ivanovich 
Kostrikin \cite{ShafKostrikin}, in which they make a very important observation 
that Cartan's classification is closely related to the classification of 
restricted Lie algebras over a field of characteristic $p> 0$. A transitive Lie algebra of a Lie pseudogroup admits a natural filtration defined by transformations that preserve $k$-jets of functions at a fixed point, which becomes an infinite-dimensional graded Lie algebra, 
or sometimes an infinite-dimensional Lie algebra. An important 
role in Cartan's classification is played by four algebras realized as 
subalgebras of the algebra of derivations of the algebra of formal power series 
$k[[t_1,\ldots,t_n]]$ over a field $k$ of characteristic 0: the algebra of all 
derivations $\mathcal{D}_n$; the algebra of all derivations $\partial$ that preserve 
the volume form $\omega =dt_1\wedge\cdots\wedge dt_n$; the algebra of 
all derivations that preserve a symplectic form; all derivations $\partial$ such 
that $\partial(\omega) = f\omega$ for some $f\in k[[t_1,\ldots,t_n]]$. 
These algebras have ideals of finite codimension that consist of 
derivations $\partial = \sum f_i\frac{\partial}{\partial t_i}$ with $f_i\in (t_1^p,\ldots,t_n)^p$. 
In characteristic $p > 0$, they represent new so-called nonclassical restricted 
Lie algebras. Kostrikin and Shafarevich made a bold conjecture 
that the class of restricted Lie algebras consists of classical ones and the 
four algebras above. In 1988, Richard Block and Robert Wilson proved this conjecture \cite{Block}. 

The study of Cartan pseudogroups led Shafarevich to investigate infinite-dimensional 
groups of biregular transformations of 
affine algebraic varieties. In his brief note \cite{ShafRendiconti} (named the ``Italian paper''), 
Shafarevich announced some fundamental results about the structure of the group of 
automorphisms of the ring of polynomials in $n$ variables based on his theory of 
infinite-dimensional algebraic groups. Answering some criticism of the lengthy review of the paper by 
Tatsuji Kambayashi, Shafarevich returns to this topic 15 years later by 
giving in \cite{ShafSome} some detailed proofs of the announced results and 
laying a foundation for the concept of an infinite-dimensional algebraic group. 
He proves that, in the case of characteristic zero, the group has a structure of a 
nonsingular infinite-dimensional algebraic variety. Another important result is 
that the group of automorphisms $\Aut(k[x_1,\ldots,x_n]$ is generated as an 
algebraic group by affine transformations and de Jonqui\'eres transformations and its subgroup 
$\Aut(k[x_1,\ldots,x_n]^0$of automorphisms with trivial jacobian is simple as 
an algebraic group. Note that neither result is true for the group of 
abstract automorphisms of the algebra. According to 
I. Shestakov and U. Umirbaev \cite{Shestakov}, the group generated by affine and 
de Jonqui\'eres transformation is a proper subgroup of $ \Aut(k[x_1,\ldots,x_3]$. and 
according to a result of Vladimir Ivanovich Danilov \cite{Danilov}, 
the group $\Aut(k[x_1,\ldots,x_2]^0$ is not simple as an abstract group. In 2004, Shafarevich 
returned to his study of infinite-dimensional groups by investigating 
the group $\textrm{GL}(2,K[t])$. He defines two different structures of an 
infinite-dimensional algebraic group on $\textrm{GL}(2,K[t])$ and studies singular 
points of their finite-dimensional closed subschemes. 

In a paper \cite{ShafDeformations}, Shafarevich studies 
the algebraic variety $\mathcal{A}_n$ parameterizing finite-dimensional 
nilpotent commutative algebras of dimension $n$ over a field. For example, 
in \cite{ShafDeformations}, he considers such algebras $N$ of nilpotent class two, 
i.e., satisfying $N^3 = 0$. In the case when the ground field is algebraically 
closed of characteristic zero, Shafarevich proves that the irreducible components of 
$\mathcal{A}_n$ coincide with its subvarieties $\mathcal{A}_{n,r}$ parameterizing algebras 
$N$ satisfying $\dim N^2= r$ assuming that $1\le r\le (n-r)(n-r+1)/2$. 
In his work, he reveals an interesting behavior of the number of irreducible 
components of $\mathcal{A}_n$.

\subsection*{Books} The name of Shafarevich is familiar to many mathematicians, 
especially to students who seek a background in algebraic geometry. 
His textbook `Basic Algebraic Geometry' was first published in Russia in 1968, 
then republished in 1972, and later published in an vastly extended version in 1988, 
and finally republished in 2007. The 1972 edition was translated into English by K.A. Hirsch in 
1974 and translated into German by Rudolf Fragel. 
The 1988 and 2007 editions were translated into English by Miles Reid in 1994 and in 2007. 

Another popular textbook written jointly with Zinovy I. Borevich is ``Theory of Numbers''. 
Its first edition was published in Russian in 1964 and republished in 1972. It was translated into 
German by Helmut Koch, into  English by Newcomb Greenleaf in  1966, and into French by 
Myriam and Jean-Luc Verley in 1967.
 
Shafarevich also published several books for a broad audience. 
A book "Geometry and groups" was written jointly with his former student Vyacheslav  Nikulin 
and published in Russian in 1983, deals with  2- and 3-dimensional locally 
Euclidean geometries and their transformation groups. 
It was translated into English by Miles Reid in 1987. 

A book ``Discourses on Algebra' translated into English by William Everett in 2003, 
is addressed to high school students and teachers. In the words of the author, 
the task of the book is to show that algebra is just as fundamental, just as deep, and just as beautiful as geometry. 

For many years, Shafarevich was one of the editors of several volumes of 
"Encyclopedia of Mathematical Sciences" published by Springer as translations 
from Russian originals published in Itogi nauki i tekhniki. Sovremennye problemy v matematike. 
Fundamentaln'ya napravleniya.
He contributed to the volumes himself by writing jointly with 
Vassily Alexeevich Iskovskikh, an article about algebraic surfaces 
in `Algebraic Geometry', vol. 3. His other contribution to the series is his book 
`Algebra I'' published in 1990 and reprinted in 1997. 
This masterpiece provides a beautiful exposition of the main concepts and 
ideas of algebra from a broader perspective of a mathematician working in various areas of mathematics. 
This confirms Shafarevich's worldview of mathematics as a unified whole, with ideas 
freely circulating from one field to another. 

\section{Non-mathematical activity}
\subsection*{Dissident movement} 
We refer to Krista Berglund's dissertation \cite{Berglund} for a meticulously 
researched, comprehensive study of this part of Shafarevich's life.  
Another rather detailed account of Shafarevich's activity as a dissident can be found in the 
book of Robert Horvath \cite{Horvath}. Here we restrict ourselves to only a brief summary of 
Shafarevich's public life outside of mathematics.

Already in 1955, Shafarevich was courageous enough to sign a letter, along with 300 other scientists, 
denouncing the works of the Soviet biologist Trofim Lysenko, who, using his power in Stalin's regime, 
opposed and prosecuted scientists working in genetics.
In 1968, Shafarevich was one of the 99 cosigners of a letter in defense of a mathematical 
logician Aleksander Esenin-Volpin, who was forcibly taken to a psychiatric hospital. 
Writing the letter deprived many of the cosigners of the possibility to travel abroad. 
Since 1971, Shafarevich has been a member of the Moscow Human Rights Committee, organized by Andrei Sakharov. 
In September 1973, he wrote an open letter in defense of Sakharov. In 1975, 
because of his dissident activity, Shafarevich was dismissed from his teaching position at 
the University (in 1949, for unknown reasons, he was also briefly dismissed from this position). 
It deprived the university of a brilliant mathematician, a popular lecturer, and a mentor of graduate students. 
As in the case of Sakharov, the membership in the Soviet Academy of Sciences and the worldwide fame
 as a scientist prevented the authorities from imposing a harsher punishment.  

In 1974, Shafarevich leaves the Sakharov Human Rights Committee and begins to collaborate with 
Alexander Solzhenitsyn in publishing an anthology   ``Is pod glyb' (`From under the rubble') \cite{Solzh1}).
First published in Russian by IMCA-Press in 1974, it was translated the following year in France, the USA, England, and Germany. 
In this collection of articles, the authors who, at that time, all resided in Russia discuss the present 
and the possible future of their country. The anthology has been condemned by the official Soviet 
propaganda as expressing the hatred of socialist ideas. The book had also been condemned by many left-leaning 
Russian dissidents as expressing Russian nationalism, chauvinism, and an attempt to replace a democratic 
society with an autocratic one. Shafarevich contributed three essays: one on ethics, one on the national problem, 
and one on socialism. The latter essay was the synopsis of his book \cite{ShafSocialism}, 
which he had already written a year before, but would publish later by the YMCA Press with a 
foreword by Solzhenitsyn in 1977. The book had been translated into French the same year. 
Earlier, before the book was released in the West, Solzhenytsyn was forcefully deported from Russia, 
So, Shafarevich had to take responsibility for discussing the book at several press conferences for 
foreign journalists (The New York Times, Frankfurter Allgemeine, BBC). 
On many occasions, Solzhenitsyn expressed his respect for Shafarevich. 
Thus, he writes in his essay ``Bodalsia telenok s dubom'' of 1975: "We have two thousand people in 
Russia, with worldwide fame, for many of them, it was much louder than for Shafarevich 
(mathematicians exist on Earth in a weak minority), however, as citizens, they are zeros 
because of their cowardice, and from this zero only a dozen took over and have grown into a tree, 
and among them is Shafarevich.'' On another occasion, he wrote: "The depth, the solidity of 
this man, not only in his figure, but in all his life image, was immediately noticed and attached.''

In 1973, Shafarevich was among a very few members of the Academy of Sciences who protested against 
the malicious campaign in the Soviet Press directed at Andrey Sakharov. He wrote an 
Open Letter distributed in Samizdat and abroad.  Next year,  he wrote two 
Open Letters protesting against the deportation of Alexander Solzhenitsyn 
with a bitter reproach to the Russian population for the unconcerned silence and even support of 
this decision. On many other occasions, Shafarevich's name could be found on various petitions in 
defense of unlawfully prosecuted human rights activists (including mathematician Leonid Plusz, 
Yuri Gastev, and physicist Yuri Osipov). Together with Sakharov, he continued to appear in court proceedings.  

After 1979, Shafarevich had stepped aside from the dissident movement. 
Although some of the dissidents tried to relate it to a crackdown on the dissident movement that started this year, 
this in no way explained by his cowardice, as his whole life amply justifies. 
As Shafarevich writes himself, he got disappointed with the movement's causes
 (like the preoccupation with the right to Jewish emigration) that he considered minor 
 compared to the real problems of the Russian people.  

\subsection*{Political activity} After Perestroika, Shafarevich began taking an active part in 
Russian political life. First supporting Yeltsyn and Sakharov, in a series of articles in "Nash Sovremennik" 
Shafarevich began to criticize the current regime for the drastic economic changes 
that left ordinary people with shortages and poverty. He also criticized 
the plans for the creation of the Soviet Sovereign Republics, which de facto should be 
dissolving the USSR. His main complaint was that this important issue needed a serious public discussion. 
The announcement of the decision had appeared five days before the date of its signature. 
The August Putsch of 1991 that followed was a tragic event (unfortunately, one of many!) in Russian history. 
In his post-putsch articles, Shafarevich compared the dissolution of the Soviet Union 
and the Communist Party with the revolution that led wide circles of ordinary people to 
despair with the new ideological and economic situation.  

As a result of this event, Shafarevich made a decision to enter politics. 
Joining the opposition camp to the regime, which was portrayed in Western media as a progressive one, 
dealt a blow to his reputation abroad. In December 1991, he joined the All-Union of Russia and 
spoke at its first congress. The new political body that united representatives of many patriotic 
and democratic movements disillusioned with Yeltsin was claimed in the West as ``the new right'', (proto)-fascist, 
and the ``red-browns''. The address of Shafarevich appealed to dropping all sectarian interests and working in the best interests of the Russian people. 
In February 1992, Shafarevich was elected (although he did not stand for election) 
to the central council of a similar new organization, the People's Gathering of Russia 
(Rossiiskoe Narodnoe Sobranie, RNS). The biased coverage of this organization by the official media, 
in particular, blaming it for the assault on its members by the Moscow TV station at Ostankino, was 
the subject of sharp criticism from Shafarevich.

In October 1992, Shafarevich joined the organizing committee of the National Salvation Front, 
representing various ideological doctrines. Very soon, by decree, Yeltsin banned the Front. 
In his statement at the Front's press conference, Shafarevich compared this with his experience as 
a dissident 20 years ago.  At that time, Yeltsin was able to consolidate his power granted to him after 
the August Putsch, and his relationship with the Congress of the Deputies (DUMA) had reached its worst. 
The statement of the organizing committee, signed by Shafarevich, demanded 
that Yeltsin and his government take responsibility for the hardship of ordinary people and suggested 
that the Front is ready to take the new executive power to prevent the country from collapsing. 

As Krista Berglund suggests 
``the moderation and sanity penetrating the Front's statement, together with lucid style and many 
formulations and emphases familiar from Shafarevich's statements,  make it plausible that he 
significantly contributed to it." The subsequent confrontation between Yeltsin and the Congress of 
the Deputies led to Yeltsin's decision to have a referendum that chose his power over the power of the Congress. 
To this referendum, Shafarevich vehemently opposed by demanding that there must be general elections for 
the President and the new Congress. As is well-known, this confrontation had ended in the bloodshed 
near the building of the Parliament 
 that left hundreds dead. Although the Front did not play any organizational role in this conflict, 
 many of its members participated in it on the side of the Parliament, compromising the Front itself. 
 After an unsuccessful attempt to be elected as a representative of the Party of Constitutional 
 Democracy in the new Parliament, Shafarevich ended his political activity. 
 Ten years later, when asked by Krista Berglund whether he had a feeling that this thing 
 [participating in political organizations] was not quite "my own", his emphatically
 agreed, except when the time he participated in the National Salvation Front.  
 After 1995, Shafarevich left all the political parties. However, since  2012, 
 he agreed to be on the editorial board of the journal ``Questions of Nationalism'' 
 of the National Democratic Party of Konstantin Krylov.

\subsection*{Non-mathematical writings} A three-volume the 
collected works of Shafarevich were published in 1994 \cite{ShafCW1}. 
In 2014, the Institute of Russian Civilization published a six-volume collected works that 
contains a lengthy introduction \cite{ShafCW2}. Only the last volume is devoted to his mathematical works. 
From the preface: "Shafarevich is a classic of Russian national thought. His books are part of the golden 
fund of Russian national heritage. For millions of Russians, the thoughts expressed in 
them become a guide in their spiritual and social life."

 Many of the non-mathematical works collected in the first five volumes were published abroad in 
 Russian or other languages.  The first such publication that appeared in the YMCA Press 
 in 1973, was the report ``Zakonodatelstvo o religii v SSSR'' (The legislation on religion in USSR'') 
 for the Human Rights  Committee. The French translation had been published in 1974 by \'Editions du Seuil, Paris. 
 His second book, "Socialism kak yavlenie mirovoy istorii" ("Socialism as a phenomenon of world history"), was published by YMCA Press in Russian in 1977 and translated 
 into French by the same publishers in the same year. Later, it was translated into 
 English as ``The socialist Phenomenon'' by Harper Collins in 1980 and published by Penguin Publications. 
 The first translation contains a preface written by A. Solzhenitsyn. 
 
 Around the same period of the seventies, Shafarevich began writing his most controversial opus "Russophobia" 
 that brought him, at the same time, love and admiration from wide circles in Russia and made 
 him a person non grata among the wide circles of Russian and Western democratic intelligentsia. 
 Although not invented by Shafarevich, the word ``Russophob'' became often associated with his book. 
 Being distributed in Samizdat in Russia since 1982, it had been officially published (in an abridged version) 
 in Russia in 1988, by a literary magazine ``Nash Sovremennik''. In the same year, the Russian 
 original was by the Munich-based journal Veche. It was 
 followed by translations into Italian (Insigna del Veltro, 1990), French (Edition Chapitre Douze, 1993), 
 Serbian (Pogledi, 1993), and  German (Verlag der Freunde, 1995). 
 It is amazing that no commercial English translation has appeared so far 
 (although Hitler's Mein Kampf is widely available both in print and on the Internet). 
 A non-commercial translation was made by Joint Publication Research Service of the 
 US Department of Commerce in 1990 and by a mathematician, Larry Shepp, in 1992 on his own initiative. 
 Never considered by Shafarevich as his most important work, the book, nevertheless, made 
 his name widely known in the West outside of the mathematical community. 
 In this book, Shafarevich borrows the theory of a French historian, Augustin Cochin (1876-1916), 
 who claimed that the French Revolution of 1789 had been initiated by a small group of intellectuals constituting 
 Malyi Narod (``Lesser or Small People'')  was opposed to the ``Large People'' who represent the organic 
 basis of the given society. Although Shafarevich did not claim 
 that "Small People" in modern Russian history consist entirely of Jews, 
 he attempted to demonstrate that the Jews indeed occupied the major part of this group. 
 As is likely to happen in any historical study, some of the factual 
 material and citations were chosen rather selectively to support his point.
 
The second volume of the collected works reprints "Russophobia" together with other 
important articles written in the nineties. Among them is one of the most important articles 
"Dve dorogi k odmomu obryvu" ("Two roads to the same abyss"). 
In this article written for the collection ``Iz pod Glyb'' which I mentioned earlier, 
Shafarevich rejects both the Socialist and the Western Democratic style for 
the future development of Russia and searches for a middle way via the spiritual reborn of the nation. 

Volume 4 of the collected works reprints another of Shafarevich's books 
"Three thousand years of mystery. History of the Jews from perspectives of modern Russia" 
published in Russia in 2002. Volume 5 contains many articles on historical and current political issues 
that appeared in the Russian Press, including three articles about Shostakovich and his music. 

Many of Shafarevich's articles were of a non-political nature, instead focusing on philosophical, 
historical, and religious topics. The leading thread of his thinking was the eternal struggle between 
good and evil. From this view, he discussed the work of Plato as well as the music of Shostakovich. 

\subsection*{Accusation in anti-semitism} The accusation is based on Shafarevich's attempt to 
defend Russia from Russophobia by expressing Judeophobia in his works. According to 
Wikipedia, anti-semitism is based on religious, economic, racist, ideological, anti-Israel, 
cultural, and social prejudices toward Jews.  Only the last one may directly apply to Shafarevich.  
The main purpose of his book, as well as of his other writings and his whole life outside mathematics, was not to 
express his hatred of Jewish people and Jewish culture, but rather to defend the Russian people, 
Russian Culture and  Russian History from accusations of their responsibility for bending under 
different political regimes, incapability to grow into a democratic society, poor cultural traditions (sic!), 
racism towards other nations and hostility to Western social ideas. 

The reaction of the mathematical community to publishing ``Russophobia'' is well known and 
widely available on the Internet. Unfortunately, the reason for the negative 
reaction of many mathematicians, many of whom probably did not bother, or were 
not able to read Shafarevich's writings,  was not the understandable concern about the 
fate of Russia in its turbulent time of the nineties,  but the outrage of what 
Shafarevich wrote concerning the Jewish people. Some of the mathematicians 
(including, for example, Jean-Pierre Serre) considered this nothing more than a witch hunt. 
Citing from a recent letter of David Mumford \cite{MumfordBlog} "I did not believe 
then and do not believe now that he was anti-semite, but rather that he was a fervent believer in his country, 
its people, its traditions -perhaps one should say its soul." For most people, the love of their country, 
its history, and its traditions, and a lesser interest or indifference to other countries and their traditions is natural. 
Unfortunately, Russia in modern times was exceptional in this way. 
The assault on the nationalistic feeling of the Russian people came from many sides: 
political, cultural, religious, intellectual, foreign, and domestic. 
Shafarevich and Solzhenitsyn were among a few people who dedicated their lives to defending the rights of 
the Russian people deserve respect from other nations.

Shafarevich expressed his own creed in the following words:
"A possibility to influence the future depends on the capability to evaluate and comprehend 
the past. Indeed, we belong to the species of Homo Sapiens, and the mind is one of 
the most powerful tools that allow us to find our own path in life. For this reason, it seems to me, 
this is now one of the most important concrete questions for Russia: stand up for the right 
to comprehend your own history without any taboo and forbidden topics.""

We may disagree with many of Shafarevich's views, some of them 
unwillingly historically distorted, but there is a good reason to remind oneself 
Voltaire's quotation: ``I disapprove of what you say, but I will defend to the death your right to say it."

Many accusations of Shafarevich being hostile to individual Jews and, especially, doing harm to 
Mathematics has not been supported by evidence. Thus, the foreign secretary of the 
Nationa Academy of Science accused Shafarevich of interfering in the careers of young 
Jewish mathematicians and preventing them from publishing their papers. 
He had never apologized for this blatant lie. 
One in four of Shafarevich's students were of Jewish, or partly Jewish, origin, and I was among them. 
Among his non-Jewish students were students of Armenian, Bulgarian, German, Litvanian, Tartar, and Ukrainian origin. 
His close associate, a friend and one of the contributors to ``Algebraic Surfaces'' was Boris Moishezon, 
one of the pioneers of the Jewish emigration movement. The coauthor of one of his most influential papers in mathematics 
 was Piatetsky-Shapiro. One of his friends (for whom he wrote a memorial article) was the 
 famous topologist Vladimir Rokhlin. 
 Shafarevich had taken a lot of effort and trouble to secure jobs for his students, 
 Jewish or not, for example, arguing before the director of the Steklov Institute, Ivan Matveich Vinogradov, 
 for the merit of giving a position at the Steklov Institute to Yuri Manin. 
 Since 1950, until his death, Shafarevich served on the editorial board of the most 
 important and prestigious Russian mathematical journal "Izvestia". 
 Between 1967 and 1977, he was the associate editor of the journal. 
 The chief editor,  Vinogradov, played only a nominal role in editorial decisions. 
 During this period, many Jewish mathematicians (e.g., Victor Kac and Boris Weisfeiler, who later 
 emigrated to the USA) were able to publish their important papers only in this journal. 

Igor Shafarevich had lived a long and productive life as a mathematician, 
a philosophical thinker, a publicist, a historian, and a Russian patriot. 
His mathematical heritage will certainly last forever; only the future will tell whether his other contributions 
to intellectual life will be of equal value.

\end{document}